\nonstopmode
\documentclass[a4paper,10pt, oneside, draft6]{amsart}
\usepackage{a4wide}
\usepackage[T1]{fontenc}
\usepackage{amsmath,mathrsfs,amssymb,amsthm,amscd,amsfonts}
\usepackage[all]{xy}
\usepackage{fancyhdr}
\usepackage{hyperref}

\newcommand{\belyi}{\boldsymbol{\beta}}
\newcommand{\NZ}{\mathbb{N}}
\newcommand{\GZ}{\mathbb{Z}}
\newcommand{\CZ}{\mathbb{C}}
\newcommand{\RZ}{\mathbb{R}}
\newcommand{\QZ}{\mathbb{Q}}
\newcommand{\OH}{\mathbb{H}}
\newcommand{\Pro}{\mathbb{P}^1}
\newcommand{\G}{\Gamma}
\newcommand{\g}{\gamma}
\newcommand{\ga}{\boldsymbol{\gamma}}
\newcommand{\si}{\sigma}
\newcommand{\SL}{SL}
\renewcommand{\div}{\operatorname{div}}
\renewcommand{\Re}{\operatorname{Re}}
\renewcommand{\Im}{\operatorname{Im}}
\DeclareMathOperator{\Stab}{Stab}

\newcommand{\mat}[4]{{\left(\begin{smallmatrix}#1&#2 \\ #3&#4\end{smallmatrix}\right)}}
\newcommand{\abcd}{{\left(\begin{smallmatrix} a&b \\ c&d\end{smallmatrix}\right)}}

\newtheorem{thm}{Theorem}[section]
\newtheorem{cor}[thm]{Corollary}
\newtheorem{prop}[thm]{Proposition}
\newtheorem{lemma}[thm]{Lemma}
\newtheorem{defin}[thm]{Definition}

\newtheoremstyle{neurem}{}{20pt}{}{}{\bf}{.}{.5em}{}
\theoremstyle{neurem}
\newtheorem{rem}[thm]{Remark}
\newenvironment{prf}{\par\pagebreak[2]\noindent{\it Proof: }}{
		      \hfill $\Box$ \medskip}
\numberwithin{equation}{thm}

\begin{document}

\markright{Posingies, \today}

\title{Kronecker limit formulas and scattering constants for Fermat curves}
\author{Anna Posingies}

\date{\today} \email{posingies@math.uni-hamburg.de}
\address{Fachbereich Mathematik (AZ)\\Universit\"at Hamburg\\
  Bundesstrasse 55\\D-20146 Hamburg\\ Germany }

\maketitle
\begin{abstract}
Eisenstein series are real analytic functions which play a central role in spectral theory of the hyperbolic Laplacian. 
Kronecker limit formulas determine their connection to modular forms.

The main result of this work is Theorem \ref{thm:grenzeinzeln} in which a Kronecker limit formula for a family of non-congruence subgroups
associated with the Fermat curves is presented.
As an application we can determine the scattering constants for the Fermat curves in Theorem \ref{thm:SCFermatall}.
\end{abstract}

\section{Introduction}

Eisenstein series are real analytic functions which play a central role in spectral theory of the hyperbolic Laplacian.
They are defined via summing over a cusp of a subgroup of $\G(1)$ (see e.g. formula \eqref{eq:ESRcuspsum}).
Kronecker limit formulas show that these functions have a strong relation to modular forms.
The classical Kronecker limit formula for $\G(1)$ is 
        \begin{align}
        \label{eq:KLFG1}
        4\pi\lim_{s\rightarrow 1}\left( E^{\G(1)}(z,s) - \frac{3/\pi}{s-1}\right)= -\log||\Delta(z)||^2  +24\left(\frac{\zeta'(-1)}{\zeta(-1)}-\log(4\pi)+1\right),
        \end{align}
(calculation similar to \cite{zagier})
where $E^{\G(1)}(z,s)$ is the Eisenstein series, $\Delta(z)$ the well known Delta function, a modular form for $\G(1)$, and $||\cdot ||^2$ the Petersson norm that will be introduced in Definition \ref{def:petnorm}.

A similar identity holds for subgroups of $\G(1)$. 
Aim of this article is to establish such a formula for the groups $\G_N$ that are associated with the Fermat curves. 
The groups $\G_N$ are of particular interest, because they are, in most cases (in all but 4), non-congruence subgroups. 
Since non-congruence subgroups are normally much harder to handle, not much is known about them.
The Fermat curves are an exceptional case due to their regularities and symmetries.
We can work with them because of their nice description (see lemmas \ref{le:grpfermat} and \ref{le:repsysGN}) and, in particular, 
because modular forms for Fermat curves were treated. 

The main result of this article is Theorem \ref{thm:grenzeinzeln} in which a Kronecker limit formula for the Fermat curves is presented.
As an application we can determine the scattering constants for the Fermat curves in Theorem \ref{thm:SCFermatall}.

\section{Eisenstein series}

\begin{lemma}
We denote by $\OH=\left\lbrace z\in \CZ \,|\, \Im(z)>0 \right\rbrace $ the upper half plane.
    The group \[\G(1):= \SL_2(\GZ)\diagup \{\pm 1\}\] acts on $\OH$ via fractional linear transformation
    \[ \abcd (z)= \frac{az+b}{cz+d} \quad \text{ for } \abcd \in \G(1) \text{ and } z\in \OH.\]
    We can add a boundary to $\OH$ by joining the upper half plane  with $\Pro(\QZ) \cong \QZ \cup \infty$, the rational projective line, and get
    $\overline{\OH}:=\OH \cup \Pro(\QZ).$
    The action of $\G(1)$ on $\OH$ can be extended to $\overline{\OH}$ by
    \[ \abcd ((p:q))= (ap+bq:cp+dq)\quad \text{ for } (p:q) \in \Pro(\QZ).\]
\end{lemma}

    \begin{prf} See \cite{miyake}.
    \end{prf}
  
At first, fix some notations. Let $\G \subset \G(1)$ be a subgroup. The classes of $\Pro(\QZ)$ with respect to the action of $\G$ are called cusps of $\G$. 
We will use the word cusp for a representative of a cusp as well.

Let $\G\subset\G(1)$ be a finite index subgroup.
For $S_j\in \Pro(\QZ)$ we will denote by $\g_j$ a matrix $\g_j\in \G(1)$ with $\g_j(\infty)=S_j$. Such a $\g_j$ always exists.
Furthermore, we normalize $\g_j$ to \[\si_j:=\g_j \cdot \mat{\sqrt{b_j}}{0}{0}{1/\sqrt{b_j}},\] with  $b_j\in \NZ$ such that 
$\sigma_j^{-1} \Gamma_j \sigma_j = \left\langle \mat{1}{1}{0}{1} \right\rangle,$
where $\G_j:=\Stab_{\G}(S_j)$.
We call $b_j$ the width of the cusp $S_j$.

For subgroups $\G\subset \G' \subset \G(1)$ a cusp $S'_k$ of $\G'$ decomposes into several cusps $\{ S_{j} \}_{j\in J_k}$ of $\G$ and
it holds $\bigcup_{j\in J_k} S_j=S'_k$. The cusps $S_j$ are the subcusps of $S'_k$ in $\G$. 
    
\begin{defin}
  \label{def:ESR}
    Let $\G \subset \G(1)$ be a finite index subgroup. 
    For each cusp $S_j$ there is a
    non-holomorphic Eisenstein series $E^{\G}_j(z,s)$, which for
    $z \in \OH$, $s \in \CZ$ and $\Re (s) >1$ is defined by the convergent series
    \begin{align}
     E^{\G}_j(z,s) = \sum_{\sigma \in \Gamma_j \setminus \Gamma} \Im
    \left( \sigma_j^{-1} \sigma(z)\right)^{s}.
    \end{align}
\end{defin}

We state some properties of Eisenstein series.
\begin{prop}
    The function $E^{\G}_j(z,s)$ has a meromorphic
    continuation to the whole $s$-plane, with a simple pole in $s=1$ of
    residue $3/( \pi \cdot [ \Gamma(1): \Gamma])=1/(vol(\G))$.\\
    Eisenstein series are automorphic forms:
    For all $\gamma\in \Gamma$ we have $E^{\G}_j(\gamma(z),s)=E^{\G}_j(z,s)$.
    They are eigenforms for the hyperbolic Laplacian $\Delta$:
    \begin{align*}
     \Delta E^{\G}_j(z,s) =s(s-1)E^{\G}_j(z,s).
    \end{align*}
\end{prop}
        \begin{prf}
        See \cite{kubota} and \cite{iwaniec}.
        \end{prf}

\begin{prop}
  \label{prop:sumESR}
    Let $\G\subset \G' \subset \G(1)$ be finite index subgroups. Let $S'_k$ be a cusp of $\G'$ and $\{ S_{j} \}_{j\in J_k}$ the subcusps of $S'_k$ in $\G$.
    The widths will be denoted by $w_k$ and $b_j$, respectively.
    Then we have the following relation for Eisenstein series
    \begin{align}
     \sum_{j\in J_k}b_j^s E^{\G}_j(z,s) = w_k^s E^{\G'}_k(z,s).
    \label{eq:sumESR}
    \end{align}
\end{prop}
        \begin{prf}
	Realizing that 
		\begin{align}
		E^{\G}_j(z,s)=b_j^{-s} \sum_{\frac{d}{c}\in S_j}\frac{\Im(z)^s}{|cz+d|^{2s}},
		\label{eq:ESRcuspsum}
		\end{align}
	the statement follows by an easy calculation. 
        \end{prf}

Eisenstein series admit a Fourier expansion. There are different normalizations that we can use in the expansion. The one below is the most useful for our purpose.

\begin{prop}
    \label{prop:expESR}  
    Eisenstein series admit a Fourier expansion.
    The Fourier expansion of $E^{\G}_j(z,s)$ at the cusp $S_k$ is given by
    \begin{align}
    E^{\G}_j(\gamma_k(z),s) &=  \delta_{jk} \frac{y^s}{b_j^s} + \pi^{1/2}\frac{\Gamma(s-1/2)}{\Gamma(s)}    \frac{1}{b_j^sb_k}\varphi_{jk,0}^{\G}(s)y^{1-s}  \notag \\
    &  \quad  + \sum_{m\neq 0} \frac{1}{b_j^sb_k} \varphi_{jk,m}^{\G}(s)   2 \pi^s \left|\frac{m}{b_k} \right|^{s-1/2} \Gamma(s)^{-1} y^{1/2} K_{s-1/2} (2\pi|m| y/b_k) e^{2\pi i m x/b_k},
    \label{eq:ESRexp}
    \end{align}
    where $z=x+iy$, $\G(\cdot)$ is the Gamma function, $K_{*}(\cdot)$ the modified Bessel function and
    \[ \varphi_{jk,m}^{\G}(s) := \sum_{c>0}\frac{1}{c^{2s}} \sum_{\substack{d \mod b_kc \\ \exists \mat{*}{*}{c}{d}\in \g_j^{-1}\G\g_k}} e^{2\pi i m \frac{d}{b_k c}}. \]
\end{prop}
	\begin{prf}
	Similar to T. Kubota \cite{kubota}. The matrix $\si_k$ in Kubota's proof has to be replaced by $\g_k$. \hphantom{box}
	\end{prf}

\begin{rem}
    \label{rem:natexpansion}
    The expansion in Proposition \ref{prop:expESR} refers to the natural cusp width by using $\g_k$ to  move the cusp. 

    A normalization using the matrix $\si_k$ (as it is done in \cite{kubota}) would lead to a different expansion, 
    which is obtained from Equation \eqref{eq:ESRexp} by replacing $z$ with $b_kz$. 
    We will call this modified expansion the normalized expansion whereas the expansion from Proposition \ref{prop:expESR} is called the natural one.
\end{rem}

We will give the following definitions using the normalized expansion. This has the advantage of coinciding with the ones in the literature
(e.g. \cite{kubota}) and of being symmetric.

\begin{defin} 
    \label{def:Cjk} For $\G\subset \G(1$) a subgroup of finite index  we define
    the scattering matrix (for the normalized Fourier expansion) to be
    \begin{align*}
    \Phi_\G(s) := \left(\pi^{1/2}
    \frac{\Gamma(s-1/2)}{\Gamma(s)} \cdot  \frac{1}{(b_jb_k)^s}\varphi^{\G}_{jk,0} \right)_{j,k},
    \end{align*}
    where $j$ and $k$ run over all cusps of $\G$.

    For all pairs $j,k$ we define the (normalized) scattering constant $C^{\G}_{jk}$ to be the constant term at
    ${s=1}$ of the Dirichlet series $(\Phi_\Gamma)_{jk}(s)$:
    \begin{align}
    C^{\G}_{jk}:= \lim_{s \to 1} \left(\Phi_\Gamma(s)_{j,k} -
    \frac{1}{vol(\G)(s-1)}   \right).
    \end{align}
\end{defin}

\begin{rem}
    \label{rem:nonnormSCvSC}
    If we take the natural Fourier expansion, see Remark \ref{rem:natexpansion}, as basis to define the (natural) scattering matrix and the (natural) scattering constants,
    they change slightly: \\
    We get $\frac{1}{b_j^s b_k}$ instead of  $\frac{1}{(b_jb_k)^s}$ in the scattering matrix. The residue does not change but the scattering constants. 
    Luckily, the difference is manageable and we have
    \begin{align}
     \widetilde{C}^{\G}_{jk} = C^{\G}_{jk}+\frac{\log(b_k)}{vol(\G)},
    \end{align}
    where $ \widetilde{C}^{\G}_{jk}$ denotes the scattering constant coming from the constant term in the natural Fourier expansion.
\end{rem}

Later on, for the Kronecker limit formulas, we will need another expansion.
From Proposition \ref{prop:expESR} follows 
\begin{cor}
	\label{cor:expgreen}
	Let $\G\subset \G(1)$ be a finite index subgroup, $S_j$ and $S_k$ cusps of $\G$.
	We have the following Fourier expansion
	\begin{multline}
	\lim_{s\rightarrow 1}\left(E^{\G}_j(\g_kz,s)- \frac{1}{vol(\G)(s-1)} \right) = \\
		\delta_{jk}\frac{y}{b_j}+\widetilde{C}_{jk}^{\G} - \frac{12\log(y)}{[\G(1):\G]}+ \sum_{m\neq 0} \frac{1}{b_jb_k} \varphi_{jk,m}(1) e^{-2\pi |m|y/b_k}\cdot e^{2\pi mx/b_k} .
	\label{eq:expgreen}
	\end{multline}
\end{cor}

For the Fourier expansions we can get a relation analogous to Equation \eqref{eq:sumESR}. Here, the calculations are more difficult and therefore we start with the preparatory 
\begin{lemma}
    \label{le:sumrs}
    Let $\G\subset \G' \subset \G(1)$ be finite index subgroups, $S_j$ represent a cusp of $\G$ as well as one of $\G'$, $S'_k$ be a cusp of $\G'$ and
    $\{ S_i \}_{i\in I_k}$ the subcusps of $\G$ such that $\cup_{i\in I_k}S_i=S'_k$. 
    By $b_{*}$ we denote the widths in $\G$ and by $w_{*}$ the ones in $\G'$, respectively.
    For $c \in \NZ$ and  an integer $m$ hold for the finite sums in the Fourier expansion of the Eisenstein series
	 \begin{align}
          \frac{1}{b_j}\sum_{l\in I_k}\sum_{\substack{d \mod b_lc \\ \exists \mat{*}{*}{c}{d}\in \g_j^{-1}\G\g_l}} e^{2\pi i m\frac{b_l}{w_k}\cdot \frac{d}{b_l c}}
	  = \frac{1}{w_j}\sum_{\substack{d \mod w_kc \\ \exists \mat{*}{*}{c}{d}\in \g_j^{-1}\G'\g_k}} e^{2\pi i m \frac{d}{w_k c}}.
	  \label{eq:sumhighcoeff}
	  \end{align}
\end{lemma} 

\begin{prf}
	In \cite{kubota} it is shown that the scattering matrix is symmetric. 
	From that, Proposition \ref{prop:sumESR} and the Fourier expansion of Eisenstein series follows that we have the desired identity for $m=0$:
	\[  \frac{1}{b_j}\sum_{l\in I_k}\sum_{\substack{d \mod b_lc \\ \exists \mat{*}{*}{c}{d}\in \g_j^{-1}\G\g_l}}1=  \frac{1}{b_j}\sum_{l\in I_k}\sum_{\substack{d \mod b_jc \\ \exists \mat{*}{*}{c}{d}\in \g_l^{-1}\G\g_j}}1 =  \frac{1}{w_j}\sum_{\substack{d \mod w_jc \\ \exists \mat{*}{*}{c}{d}\in \g_k^{-1}\G'\g_j}}1 = \frac{1}{w_j}\sum_{\substack{d \mod w_kc \\ \exists \mat{*}{*}{c}{d}\in \g_j^{-1}\G'\g_k}}1\]
	That gives us the number of summands in the double sum on the left in Equation \eqref{eq:sumhighcoeff} comparative to the sum on the right:
	$(\text{sums left}) = \frac{b_j}{w_j}\cdot (\text{sum right})$.

	There is a well known decomposition of $\g_j^{-1}\G\g_k$ into double cosets (see \cite{iwaniecclass} slightly modified):
	\begin{align}
	\g_j^{-1}\G\g_k = \delta_{jk}\left \langle\mat{1}{b_j}{0}{1}\right \rangle  \cup \bigcup_{c\ge 0} \bigcup_{d \text{ mod } b_kc} \left\langle\mat{1}{b_j}{0}{1}\right \rangle \mat{*}{*}{c}{d}\left \langle\mat{1}{b_k}{0}{1}\right \rangle,
	\end{align}
	where the union is taken over all pairs $c,d$ such that there is $\mat{*}{*}{c}{d} \in \g_j^{-1}\G\g_k$.
	Since each $d$ in the decomposition determines the upper left entry of the matrix $a \mod b_jc$, we get   
	\begin{align}
	\g_j^{-1}\G\g_k = \delta_{jk}\left \langle\mat{1}{b_j}{0}{1}\right \rangle  \cup \bigcup_{c\ge 0} \bigcup_{a \text{ mod } b_jc} \left\langle\mat{1}{b_j}{0}{1}\right \rangle \mat{a}{*}{c}{*}\left \langle\mat{1}{b_k}{0}{1}\right \rangle. \label{eq:decompas}
	\end{align}
	Hence, the $d'$s on the left hand side of Equation \eqref{eq:sumhighcoeff} correspond to $a's$ in matrices $\mat{a}{*}{c}{d}$. 
	With the decomposition \eqref{eq:decompas}, one can show that all $a$'s in the matrices are distinct. 

	From $\mat{a}{b}{c}{d}\in \g_l^{-1}\G\g_j$ follows $\mat{a}{b}{c}{d}\in \g_k^{-1}\G'\g_j$. 
	Thus, the reduction of $a\mod w_j$ creates a class in the double coset decomposition of $\G'$.  
        Hence, all $a$'s lie in
	\[ A =  \left \lbrace a'+nw_j \,\big |\,  0\le a' < w_jc \, \text{ such that }\,\exists \mat{a'}{*}{c}{*} \in \g_j^{-1}\G'\g_k, 0\le n < \frac{b_j}{w_j}\right \rbrace.\]  
	Since $|A|$ is just the required number, the whole of $A$ is the set of $a$'s. Now we see that the reduction modulo $w_j$ of $\frac{b_j}{w_j}$ different $a$'s coincide, the reduction of the corresponding $d$'s likewise. Hence, the statement follows.
\end{prf}

Now, we can prove
\begin{prop}
    \label{prop:sumexpgreen}
    Let $\G \subset \G' \subset \G(1)$ be finite index subgroups. Let $S_j$ and $S_k$ be cusps of cusp width $b_*$ (in $\G$) and $w_*$ (in $\G'$), respectively.
    Then we have
    \begin{align}
       \frac{1}{b_j^{1-s}}\sum_{\g\in\G\setminus\G'} E_j^{\G}(\g\g_k z,s) =\frac{1}{w_j^{1-s}}E_j^{\G'}(\g_k z,s). 
    \end{align}
\end{prop}
        \begin{prf}
	To understand what happens in the sum, we choose a suitable system of representatives.
        We have
        \[ \G' = \bigcup_{S_l \sim_{\G'} S_k } \bigcup_{0\le n < \frac{b_l}{w_k}} \G \g_{lk}\tau_{k,w_kn},\]
	where the first union is taken over representatives for all cusps $S_l$ of $\G$ that are $\G'$-equivalent to $S_k$.
        The matrices that occur are  all from $\G(1)$ with  $\g_{lk}(S_k)=S_l$,  $\tau_{k,w_kn}=\g_k\tau_{w_kn}\g_k^{-1} \in  \Stab_{\G'}(S_k)$ and $\tau_{w_kn}=\mat{1}{w_kn}{0}{1}$.
        Then
        \begin{align*}
        \sum_{\g\in\G\setminus\G'}E_j^{\G}(\g\g_k z,s)
        &=\sum_{S_l \sim S_k } \sum_{0\le n < \frac{b_l}{w_k}}
        E_j^{\G}(\g_{lk}\g_k\tau_{w_kn}\g_k^{-1}\g_k z,s) \\
        &=\sum_{S_l \sim S_k } \sum_{0\le n < \frac{b_l}{w_k}}  E_j^{\G}(\g_{l}(z+w_kn),s)
        \end{align*}
        with $\g_{lk}\g_k=\g_l$ which fulfills $\g_l(\infty)=S_l$. 

   	Now, we look at the sum of the Fourier expansions (use Equation \eqref{eq:ESRexp}) and we get 
	\begin{align*}
	\sum_{\g\in\G\setminus\G'}E_j^{\G}(\g\g_k z,s) &= \delta_{jk}\frac{b_j}{b_j^sw_j}y^s+\pi^{1/2}\frac{\G(s-1/2)}{\G(s)}\frac{b_j}{b_j^sw_jw_k} \left(\sum_{c> 0}\frac{1}{c^{2s}}\sum_{S_l\sim S_k} \sum_{\substack{d \mod b_lc \\ \exists \mat{*}{*}{c}{d}\in \g_j^{-1}\G\g_l}} 1\right)y^{1-s} \\
	&\quad + \sum_{m\neq 0} \sum_{S_l\sim S_k} \frac{1}{b_j^sb_l}\left(\sum_{c>0}\frac{1}{c^{2s}} \sum_{\substack{d \mod b_lc \\ \exists \mat{*}{*}{c}{d}\in \g_j^{-1}\G\g_l}} e^{2\pi i m \frac{d}{b_l c}}\right)  2 \pi^s \left|\frac{m}{b_l} \right|^{s-1/2} \\
	&\quad\cdot  \Gamma(s)^{-1} y^{1/2} K_{s-1/2} (2\pi|m| y/b_l) \left(\sum_{0\le n < \frac{b_l}{w_k}}e^{2\pi i n /(b_l/w_k)} \right) e^{2\pi i m x/b_l}.
	\end{align*}
	On the constant term we can apply Lemma \ref{le:sumrs}. For the higher terms we first conclude with
        \[ \sum_{0\le n < k}e^{2\pi i m \frac{n}{k}} = 
				 \begin{cases}
                                 k & \text{ if } k|m \\
                                 0 & \text{ elsewise,}
                                 \end{cases}   \]
	that for many $m$ the coefficient is zero and apply Lemma \ref{le:sumrs} on the remaining ones. 
	If we compare the result with the Fourier expansion of $E_j^{\G'}(\g_kz,s)$ we get the statement.  
	\end{prf}

Similarly, an identity holds for the values in $s=1$.

\begin{cor}
   \label{cor:sumexpgreen}
   With the notations from Proposition \ref{prop:sumexpgreen} we have
   \begin{multline}
       \sum_{\g\in\G\setminus\G'}\lim_{s\rightarrow 1}\left(E_j^{\G}(\g\g_k z,s)-\frac{1}{vol(\G)(s-1)}\right) =\lim_{s\rightarrow 1}\left(E_j^{\G'}(\g_k z,s)-\frac{1}{vol(\G')(s-1)}\right) \\ -\frac{1}{vol(\G')}\cdot \log\left( \frac{b_j}{w_j}\right).
   \end{multline}
\end{cor}

\section{Kronecker limit formulas}

To establish Kronecker limit formulas we need modular forms. Here, we just give the definitions needed. An introduction to modular forms can be found in \cite{miyake}.

\begin{defin}
    \label{def:slashop}
    We define an action of $\G(1)$  on functions $f:\OH \rightarrow \CZ$ via
    \[ f|_k \g (z)= (cz+d)^{-k}f(\g z), \]
    where $\g=\abcd\in\G(1)$ and $k\in \GZ$. This action is called the slash operator of weight $k$ or the $k$-th slash operator.
\end{defin}

Let $\G\subset \G(1)$ be a subgroup of finite index, $k$ an integer. 
A meromorphic function $f: \OH \rightarrow \CZ$ behaves automorphically of weight $k$ with respect to $\G$  if
    \[ f|_k \g (z)= f(z) \qquad \forall \g\in \G.\]
Then $f|_k\g_j (z)$ is $b_j$-periodic, i.e. $ f|_k\g_j (z+b_j)=f|_k\g_j (z)$, where $b_j$ is the cusp width of $S_j$. 
Therefore there exists a function $g$ on $D\setminus \{0\}$ (the punctured unit disc) such that
    $$ f|_k \g_j (z) =g(e^{2\pi i z}) \qquad z\in \OH. $$
The function $g$ is meromorphic on $D\setminus \{0\}$, since $f$ is meromorphic. 
We say that $f$ is meromorphic, is holomorphic in the cusp $S_j$ if $g$ extends meromorphically, holomorphically to $0$, respectively.

\begin{defin}
    \label{def:modformfunction}
    A holomorphic function $f: \OH \rightarrow \CZ$ is called modular function with respect to $\G$ if it behaves automorphically of weight 0 and is meromorphic in all cusps of $\G$;
    a holomorphic function  $f(z)$ is called a modular form (of weight $k$ with respect to $\G$) if it behaves automorphically of weight $k$ and is holomorphic in all cusps.

    The set of modular forms of weight $k$ with respect to $\G$ is denoted by  $M_k(\G)$, it generates a ring graded by the weight.
\end{defin}
 
\begin{defin}
        \label{def:petnorm}
        Let $f(z)\in M_k(\G)$ be a modular form for a finite index subgroup $\G\subset \G(1)$. Then we define its Petersson norm via
	\begin{align}
        || f(z)||^2:= |f(z)|^2 \Im(z)^k.
        \end{align}
\end{defin}

Now, we can formulate the Kronecker limit formula:
\begin{prop}
    \label{prop:KLF}
    Let  $\G\subset \G(1)$ be a subgroup, $S_j$ a cusp of $\G$. Suppose there is a modular form $f_j^{\G}\in M_{k}(\G)$ ($k\in \NZ$), that only vanishes in the cusp $S_j$.
    Then there is a constant $A\in \RZ$ (depending on $f_j^{\G}$) such that
   \begin{align}
        \label{eq:KLF}
     4\pi \lim_{s\rightarrow 1}\left(E^{\G}_j(z,s)-\frac{1}{vol(\G)(s-1)}\right) = - \frac{1}{[\G(1):\G]\cdot\frac{k}{12}}\log || f_j^{\G} ||^2 +A.
    \end{align}
\end{prop}
	\begin{prf}
        Examine the action of the hyperbolic Laplace operator $\Delta=y^2\left(\frac{\partial^2}{\partial x^2}+\frac{\partial^2}{\partial y^2} \right)$ on the functions $4\pi \lim_{s\rightarrow 1}\left(E^{\G}_j(z,s)-\frac{1}{vol(\G)(s-1)}\right)$ and  $-\log || f_j^{\G}||^2$. 
	We have for the left side of Equation \eqref{eq:KLF} (for the expansion in a cusp $S_l$)
        \begin{align*}
        &\Delta\left(4\pi \lim_{s\rightarrow 1}\left(E^{\G}_j(\g_l z,s)-\frac{1}{vol(\G)(s-1)}\right)\right)\\
	&\hspace{3cm}= \Delta \left(4 \pi \left( \delta_{jl}\frac{y}{b_l} + \widetilde{C}^{\G}_{jl}-\frac{3\log(y)}{\pi [\G(1):\G]}+\sum_{n\neq 0} a_n(y,1)q^{xn/b_l} \right) \right) \\
	&\hspace{3cm}= \frac{12}{ [\Gamma(1):\Gamma] },
        \end{align*}
	with $q=e^{2\pi i}$ and $z=x+iy$, where the formula for the expansion can be found in Corollary \ref{cor:expgreen}.
        On the other side,
        \begin{align*}
        \Delta(-\log|| f_j^{\G}||^2) &=\Delta \left(-\log \left(|f_j^{\G}|^2 y^{k}\right) \right) \\
        &=   -\Delta \left(\log(f_j^{\G}) \right)-\Delta(\log(\overline{f_j^{\G}}))- \Delta \left(\log \left(y^{k}\right)\right)\\
        &=  0+0+ k.
        \end{align*}
        Hence: \[\Delta\left(4\pi \lim_{s\rightarrow 1}\left(E^{\G}_j(z,s)-\frac{1}{vol(\G)(s-1)}\right)+\frac{1}{[\G(1):\G]\cdot \frac{k}{12}} \log || f_j^{\G} ||^2\right)=0 \] 
		The spectral decomposition of the Laplacian had been studied, see \cite{iwaniec}, for functions that are square integrable (the space $\mathfrak{L}(Y_{\G})$). 
		To find out if we can use the result from \cite{iwaniec}, we will study the behavior of $ 4\pi \lim_{s\rightarrow 1}\left(E^{\G}_j(z,s)-\frac{1}{vol(\G)(s-1)}\right)+\frac{1}{[\G(1):\G]\cdot \frac{k}{12}} \log || f_j^{\G} ||^2 $ in the cusps.
		For that, we will compare the expansions. We have seen the expansion of the Eisenstein series in a cusp $S_l$ in Corollary \ref{cor:expgreen}. 
	The expansion $f_j^{\G}|_{S_l}$ has the form $d_m q^{zm/b_l} \left ( 1+\sum_n d_nq^{zn/b_l} \right)$, with $q=e^{2\pi i }$,
	$m\in \left\{0,[\G(1):\G]\cdot \frac{k}{12}\right\}$ and $m=0$ if and only if $S_j\not=S_l$, since $f_j^{\G}$ only vanishes in the cusp $S_j$ and the vanishing order is $[\G(1):\G]\cdot \frac{k}{12}$ (by the theory of modular forms).
	Therefore, we have ($z=x+iy$)
	{\allowdisplaybreaks	
	\begin{align*}
        \log || f_j|_{S_l}||^2 &= \log \left( y^{k}| f_j|_{S_l}|^2 \right) \vphantom{\sum_n}\\
            &= k\cdot\log(y)+ 2\Re \log (f_j|_{S_l}) \\
            &=  k\cdot\log(y)+ 2\Re \log \left(d_mq^{zm/b_l} \left ( 1+\sum_{n>0} d_nq^{zn/b_l} \right)\right) \\
            &=  k\cdot\log(y)+ 2\Re \log(d_m)-\delta_{jl}4\pi \frac{y}{b_l}\cdot [\G(1):\G]\cdot \frac{k}{12}\\ 
		& \quad +2\Re \log \left ( 1+\sum_{n>0} d_nq^{zn/b_l} \right) \\
            &= k\cdot\log(y)+ 2\Re \log(d_m)-\delta_{jl}4\pi \frac{y}{b_l}\cdot [\G(1):\G]\cdot \frac{k}{12}
			 +2\Re \left ( \sum_{n>0} \tilde d_nq^{zn/b_l} \right),
        \end{align*} }
	 where the $\tilde d_n$ are suitable such that $\log \left ( 1+\sum_{n>0} d_nq^{zn/b_l} \right)=\sum_{n>0} \tilde d_nq^{zn/b_l}$.

	Now we can see that the value of \[ 4\pi \lim_{s\rightarrow 1}\left(E^{\G}_j(z,s)-\frac{1}{vol(\G)(s-1)}\right)+ \frac{1}{[\G(1):\G]\cdot \frac{k}{12}}\log || f_j^{\G} ||^2 \] in $S_l$ is bounded:
	\begin{align*}
        & \lim_{z\rightarrow i\infty}\left(4\pi \lim_{s\rightarrow 1}\left(E^{\G}_j(\g_lz,s)-\frac{1}{vol(\G)(s-1)}\right)+ \frac{1}{[\G(1):\G]\cdot \frac{k}{12}}\log || f_j^{\G}|_{S_l} ||^2\right) \\
        &\quad \qquad= \lim_{z\rightarrow i\infty}\left(4\pi \widetilde{C}_{jl}^{\G}+\frac{12\log(4\pi)}{[\G(1):\G]}+2\Re \log(d_m)+ \sum_{n\neq 0} a_n e^{-2\pi |n|\frac{y}{b_l}}e^{2\pi i \frac{x}{b_l}} \right. \\
	&\quad \qquad\quad \left.+2\Re \left ( \sum_{n>0} \tilde d_ne^{2\pi in\frac{z}{b_l}} \right) \right) \\
        &\quad \qquad  = 4\pi \widetilde{C}_{jl}^{\G}+\frac{12\log(4\pi)}{[\G(1):\G]}+2\Re \log(d_m),
        \end{align*}
        where $S_l$ was chosen arbitrarily. Therefore, we have  \[4\pi \lim_{s\rightarrow 1}\left(E^{\G}_j(z,s)-\frac{1}{vol(\G)(s-1)}\right)+\frac{1}{[\G(1):\G]\cdot \frac{k}{12}} \log || f_j ||^2 \in \mathfrak{L}(Y_{\G}).\]
        The spectral decomposition in  \cite{iwaniec} shows that the kernel of the Laplacian are the constant functions. 
        If we take a closer look at the functions involved here, we see that we get a real number. 
\end{prf}

\begin{rem}
   \label{rem:Aindepcusp}
    The constant $A$ from Proposition \ref{prop:KLF} in Formula \eqref{eq:KLF} can be calculated by comparison of the Fourier expansions. 

    It is independent of the cusp $S_l$ in which the expansion is taken:
    To explain this, we may start with the functions expanded in $\infty$ to see what happens if we pass on to another cusp.  Changing from $\infty$ to the cusp $S_l$ means changing $z$ to $\g_lz$, where $\g_l\in \G(1)$ with $\g_l^{-1} \G_l \g_l=\left \langle \mat{1}{b_l}{0}{1} \right \rangle $.
    Thus, only the parts of the function change that depend on $z$. These parts coincide for both sides of Equation \ref{eq:KLF} such that their difference stays the same.
\end{rem}

\section{Kronecker limit formulas for $\G(2)$}

The group $\G(2)\subset \G(1)$ is a free subgroup of index $6$ with two generators
	\begin{align} \ga_1:=\mat{1}{2}{0}{1} \quad \text{and} \quad  \ga_2:=\mat{1}{0}{2}{1}. \label{eq:gensG2}\end{align}
It has three cusps of width $2$ that are (see \cite{shimura}):
    \begin{align}
    &\left\lbrace (p:q)\in \Pro(\QZ) | (p,q)=1, (p,q)\equiv (0,1) \mod 2 \right\rbrace \notag \\
    &\left\lbrace (p:q)\in \Pro(\QZ) | (p,q)=1, (p,q)\equiv (1,1) \mod 2  \right\rbrace  \label{eq:cuspsG2} \\
    &\left\lbrace (p:q)\in \Pro(\QZ) | (p,q)=1, (p,q)\equiv (1,0) \mod 2  \right\rbrace, \notag
    \end{align}
    thus a system of representatives is $\{ 0, 1, \infty \}$.

 We define (with $q=e^{2\pi i z}$):
    \begin{align}
    \theta^2(z) &= \prod_{n\ge 1}\left( 1-q^n\right)^4\left( 1+q^{n-1/2}\right)^8  \notag \\
    \lambda(z) &= -\frac{1}{16} q^{-1/2} \prod_{n\ge 1}\left(\frac{1-q^{n-1/2}}{1+q^n} \right)^8 \label{eq:prodmodforms}\\
    (1-\lambda)(z) &= \frac{1}{16} q^{-1/2} \prod_{n\ge 1}\left(\frac{1+q^{n-1/2}}{1+q^n} \right)^8 \notag
    \end{align}
    Then $\theta^2(z)$ is a modular form for $\G(2)$ of weight $2$, the functions $\lambda(z)$ and $(1-\lambda)(z)$ are modular functions for the same group.
    They have the following divisors
    \[ \div \theta^2 = 1\cdot 1, \quad \div \lambda = 1\cdot 0 - 1\cdot \infty , \quad \div (1- \lambda) = 1\cdot 1 - 1\cdot \infty.\]
    (See \cite{yang}, or, for more background information, \cite{miyake} and \cite{emot}.)

Hence, the modular forms
        \begin{align}
        G_0(z) &:= \frac{\lambda(z)}{1-\lambda(z)}\theta^2(z) \label{eq:mofformG0}\\
        G_1(z) &:= \theta^2(z) \label{eq:mofformG1}\\
        G_{\infty}(z) &:= \frac{1}{1-\lambda(z)}\theta^2(z) \label{eq:mofformGinf}
        \end{align}
        from $M_2(\G(2))$ have divisors
        \[ \div G_0 = 1\cdot 0, \quad \div G_1 = 1\cdot 1 , \quad \div G_{\infty} = 1\cdot \infty.\]

\begin{prop}[Kronecker limit formula for $\G(2)$]
	\label{prop:KLFG2}
	For the group $\G(2)$ holds
	\begin{multline}
        4\pi \lim_{s\rightarrow 1}\left( E^{\G(2)}_j(z,s)-\frac{1}{vol(\G(2)) (s-1)}\right)= \\ -\log|| G_j(z)||^2+
        4\left(\frac{\zeta'(-1)}{\zeta(-1)}-\log(4\pi)+1+\frac{1}{6}\log(2)\right),
	\label{eq:KLFG2}
        \end{multline}
        where $j\in \{0,1, \infty\}$ denotes one of the three cusps of $\G(2)$, $E^{\G(2)}_j(z,s)$ is an Eisenstein series and
        $G_j$ the corresponding modular form from one of the equations \eqref{eq:mofformG0} to \eqref{eq:mofformGinf}.
\end{prop}
    \begin{prf}
    We will compare the expansions of both functions involved in the cusp $\infty$.
    In $\infty$  we get expansions (with $z=x+iy$)
    \[4\pi \lim_{s\rightarrow 1}\left( E^{\G(2)}_j(z,s)-\frac{1}{vol(\G(2)) (s-1)}\right) = \sum_{m\in \GZ}e_{j,m}(y)e^{\pi i m x} \]
    and
    \[-\log|| G_j(z)||^2 = \sum_{m\in \GZ}g_{j,m}(y)e^{\pi i m x}. \]
    The identity $e_{j,m}(y)=g_{j,m}(y)$ for all $m\neq 0$ follows from Proposition \ref{prop:KLF}.
    Therefore, we just have to deal with $m=0$.
    The coefficient $g_{j,0}(y)$ can easily be derived from the product description in Equations \eqref{eq:prodmodforms} 
    and the definition of the Petersson norm (\ref{def:petnorm}). In all three cusps $S_j$ holds
    \begin{align} g_{j,0}(y)=\delta_{j\infty}\left(2\pi y-8\log(2)\right)-2\log(y). \label{eq:consttermfform}\end{align}
    If we regard the Eisenstein series, we get from Equation \eqref{eq:expgreen} that 
    \begin{align} e_{j,0}(y)=\delta_{j\infty}2\pi y +4\pi \widetilde{C}_{j\infty}^{\G(2)}-2\log(y). \label{eq:consttermfesr}\end{align}

    The scattering constants for the group $\G(2)$ can be calculated with the description of the cusps (Equation \eqref{eq:cuspsG2}), either by using results by Huxley \cite{huxley} or directly from the Fourier expansion (as it had been done in \cite{posingiesdiss}): One gets
    \begin{align*}
      \widetilde{C}^{\G(2)}_{jk}& =\lim_{s\rightarrow 1} \left( \frac{1}{2\cdot 2^s}\pi^{1/2}\frac{\G(s-1/2)}{\G(s)}\frac{2^{2s}-2}{2^{2s}-1}\frac{\zeta(2s-1)}{\zeta(2s)}-\frac{1}{2\pi (s-1)}\right) \\
      &=\frac{1}{\pi}\left(\frac{\zeta'(-1)}{\zeta(-1)} -\log(4\pi)+1+\frac{1}{6}\log(2) \right) \\
      \widetilde{C}^{\G(2)}_{jj} &=\lim_{s\rightarrow 1} \left( \frac{1}{ 2^s}\pi^{1/2}\frac{\G(s-1/2)}{\G(s)}\frac{1}{2^{2s}-1}\frac{\zeta(2s-1)}{\zeta(2s)}-\frac{1}{2\pi (s-1)}\right)\\
      &=\frac{1}{\pi}\left(\frac{\zeta'(-1)}{\zeta(-1)} -\log(4\pi)+1-\frac{11}{6}\log(2) \right),
    \end{align*}
    where the first formula holds when $S_j\neq S_k$ and the second one in case of equality. 
    With this information we can compare equations \eqref{eq:consttermfform} and \eqref{eq:consttermfesr} to obtain the statement.
\end{prf}

\section{Basics on Fermat curves}

As a projective curve the well known Fermat curve is given by
\begin{defin}
    Let $N \in \NZ$. The $N$-th Fermat curve is given by the equation
    \begin{align} \label{eq:fermatcurve} F_N: \quad X^N+Y^N = Z^N. \end{align}
\end{defin}

\begin{lemma}
    \label{le:BelyiFermat}
    Consider the map
    \begin{align} \label{eq:fermatbelyi}
     \belyi_N: \quad  F_N &\longrightarrow \Pro. \\
        (X:Y:Z) & \longmapsto (X^N:Z^N) \notag
    \end{align}
    Its degree is $N^2$. It is ramified only above the points $0,1,\infty$ and the ramification points are
    \begin{align}
     a_j &:= (0:\zeta^j: 1) \notag \\
     b_j &:= (\zeta^j:0: 1) \label{eq:cuspsfermat}\\
     c_j &:= (\epsilon \zeta^j:1: 0), \notag
    \end{align}
    where $\zeta=e^{2\pi i /N}$ is the first primitive $N$-th root of unity, $ j\in \{ 0, \dots N-1\}$ and $\epsilon=e^{\pi i/N}$. Each point has ramification index $N$.
\end{lemma}
    \begin{prf}
        Simple calculation.
    \end{prf}

There is a subgroup $\G_N$  of $\G(2)$ given by the monodromy of the cover $\belyi_N$  with the property \[ \G_N \setminus \OH \cong F_N(\CZ)\setminus\{\text{ramification points of }\belyi_N\}.\]
The group $\G_N$ can be described as in the following

\begin{lemma}
    \label{le:grpfermat}
    The group  $\G_N$ is the kernel of
    \begin{align*}
    \G(2) &\longrightarrow \GZ / N\GZ \times \GZ / N\GZ, \\
    \g & \longmapsto (R_1(\g), R_2(\g)) \mod N,
    \end{align*}
    where $R_i(\g)$ denotes the number of generators $\ga_i$ ($i\in \{0,1\}$) for $\G(2)$ (see Equation \eqref{eq:gensG2}) that occur in the word description of $\g$: \\
     Let $\g\in\G(2)$ be given via its word in $\ga_1$ and $ \ga_2$ as
     $\g=\prod_{i=1}^n \kappa_i^{r_i}$ with $n\in \NZ, \; r_i\in \GZ$,  \mbox{$\kappa_i\in \{\ga_1, \ga_2\}$.} Then
    \[ R_1(\g)= \sum_{\kappa_i=\ga_1}r_i \quad \text{and} \quad R_2(\g)= \sum_{\kappa_i=\ga_2}r_i.\]
\end{lemma}
    \begin{prf}
        See \cite{murtyrama}.
    \end{prf}

\begin{rem}
     The subgroup $\G_N$ from Lemma \ref{le:grpfermat} is a non-congruence subgroup for all $N$ but $1,2,4$ and $8$ (see \cite{phillipssarnak}). 
\end{rem}

Further facts about $\G_{N}$.
\begin{lemma}
        \label{le:repsysGN}
        We have $\G_N\vartriangleleft \G(1)$, $\left[\Gamma(1):\Gamma_{N} \right] = 6N^2$ and the group has $3N$ cusps, all of same width $b=2N$.
        A system of representatives for the cosets $\G_N \setminus \G(2)$ is
        \begin{align}
        \{\ga_1^a \ga_2^b\} \qquad \text{with} \qquad a,b\in {0, \dots , N-1}.
        \label{eq:cosetsG2}
        \end{align}
 	A system of representatives for the cusps is $S=S_0 \cup S_1 \cup S_{\infty}$ with
        \begin{align}
        \label{eq:repsysFN}
        S_0=\left\lbrace 0,2,\dots, 2N-2\right\rbrace, S_1=\left\lbrace  1, 3, \dots, 2N-1\right\rbrace, S_{\infty}=\left\lbrace  \frac{1}{2}, \frac{1}{4}, \dots, \frac{1}{2N} \sim_{_{\G_N}} \infty\right\rbrace . \end{align}
        The cusps in $S_i$ are $\G(2)$-equivalent to $i$ ($i\in\{0,1, \infty\}$).
\end{lemma}
	\begin{prf}
	The fact that $\G_N$ is normal, the index, the number of cusps and the representatives for the cosets follow from the lemmas \ref{le:BelyiFermat} and \ref{le:grpfermat}.	

	To prove that $S$ contains exactly one representative for all cusps, it is enough to show that all elements of $S$ are non-equivalent under $\G_N$. Since cusps from different subsets $S_0, S_1, S_{\infty}$ are non-equivalent under $\G(2)$ (see Equation \eqref{eq:cuspsG2}), we have to compare cusps out of the same subset only. 
	Hence, we have to express a general matrix that maps cusps out of one subset to each other in the generators $\ga_1, \ga_2$   of $\G(2)$ and find out if they are in $\G_N$. 

	This reduces the problem to the examination of the following words, where $j$ and $k$ denote cusps from $S_i$, $i\in \{0,1, \infty\}$: 
	\[ \g_{ik}^{-1}\kappa_i^m\g_{ij} \quad \text{with} \quad m\in\GZ;\, \Stab_{\G(2)}(i)=\left\langle \kappa_i \right \rangle; \g_{ij},\g_{ik}\in\G(2):  \g_{ij}(j)=\g_{ik}(k)=i. \] 
	We have $\kappa_0=\ga_2, \kappa_1=\big(\ga_2\ga_1^{-1}\big),\kappa_{\infty}=\ga_1$. In the cases $i=0$ and $i=1$ powers of $\ga_1$ do for $\g_{ki}$ as well as for $\g_{ji}$ and for $i=\infty$ powers of $\ga_2$. 
	Combined with the fact that the smallest value of $|m|$, for which $\ga_1^m$ or $\ga_2^m$ lie in $\G_N$, is $|m|=N$, we get the result.	 
	\end{prf}

There is a 1-1-correspondence between the cusps of $\G_N$ and the ramification points of $\belyi_N$ that can be made explicit.

\begin{prop}
	\label{prop:corrcuspsfermat}
	A possible identification of ramification points of the Belyi map $\belyi_N$  (Equation \eqref{eq:cuspsfermat}) and the cusps of the group $\G_N$ (Equation \eqref{eq:repsysFN}) is
	\[\begin{array}{cccc}
	   (0:\zeta^n: 1)  & \longleftrightarrow & 2N-2n &\\ 
	(\zeta^n: 0 : 1)  & \longleftrightarrow & 2N-2n-1 &\quad n\in\{0, \dots ,N-1\}. \\ 
        (\epsilon\zeta^n: 1: 0) & \longleftrightarrow & \frac{1}{2N-2n} &
	\end{array}\]
\end{prop}
	\begin{prf}
	The situation is the following:
	\begin{align*}
        \begin{xy}
        \xymatrix{
        F_N(\CZ)\setminus \{a_j, b_j, c_j\}_{j=0\dots N-1} \ar@{=}[r]^/1.8em/{\sim}\ar[d] & \G_N\setminus \OH \ar[d]\\
       \Pro(\CZ)\setminus\{0,1,\infty\} \ar@{=}[r]^/.6em/{\sim} & \G(2)\setminus \OH
        }
        \end{xy}
        \end{align*}
	The lower isomorphism is given via $\lambda(x)$ (see equation \eqref{eq:prodmodforms}; we compose it with an automorphism to fix $0,1,\infty$). 
	The isomorphism preserves the orientation.

	Via lifting $N$-fold circles around $0$ and $\infty$  we get the orders of the cusps around $(0:1: 1)$ and $(\epsilon: 1: 0)$.

	For $(0:1: 1)$ take the path $\frac{1}{2}e^{2\pi i N \lambda}$ ($\lambda\in [0,1]$) on $\Pro$. 
        One of its lifts on $F_N$ under $\belyi_N$ meets all the preimages of the real line that connect the cusp $(0:1: 1)$ with the ${b_j}'s$ and ${c_j}'s$ 
	(Equation \eqref{eq:cuspsfermat}) in the following order:
	\[ (1:0:1), (\epsilon: 1 :0), (\zeta_N:0:1), (\epsilon \zeta_N:1:0), \dots , (\zeta_N^{N-1}:0:1), (\epsilon \zeta_N^{N-1}:1:0) .\]
	For $\infty$ we can get a corresponding result by lifting $2e^{2\pi i N \lambda}$: 
	\[(1:0:1), (0:1:1),(\zeta_N:0:1) ,  (0:\zeta_N:1), \dots ,  (\zeta_N^{N-1}:0:1) ,(0:\zeta_N^{N-1}:1) .\]
 
	A lift of a circle around the cusp $S_i$ on the $F_N$ side corresponds to the application of a matrix $\kappa_i$, 
	that generates $\Stab_{\G(2)}(S_i)$, on the side of $ \G_N\setminus \OH$. The quotient $ \G_N\setminus \OH$ is represented by a fundamental domain $\mathcal F_N\subset \OH$ for $\G_N$. 
	It has a tessellation by fundamental domains of $\G(2)$, that are triangles with vertices $0$, $1$ and $\infty$. 
	When $\kappa_i$ acts on $\mathcal F_N$ it interchanges the triangles around the cusp $S_i$ and we get an order of the cusps as before.

	In the case of $0$ we have to take $\kappa_0=\ga_2^{-1}$ to describe a turn in positive direction. 
	Then we get $\ga_2^{-n}(\infty)=\frac{1}{-2n}\sim_{_{\G_N}}\frac{1}{2N-2n}$ (the equivalence is given via $\ga_2^N\in \G_N \quad \forall N$)
	as well as $\ga_2^{-n}(1)=\frac{1}{-2n+1}\sim_{_{\G_N}}2N-2n+1$ (via $\ga_1^{N+n-1}(\ga_2\ga_1^{-1})^{n-1}\ga_2^{1-n}\in \G_N \quad \forall N$).

	Therefore the order of the cusps around $0$ is:
	\[ \ga_2^0(1), \ga_2^0(\infty), \ga_2^{-1}(1), \ga_2^{-1}(\infty),  \dots  , \ga_2^{-N+1}(1) ,\ga_2^{-N+1}(\infty) =  1,\infty, 2N-1, \frac{1}{2N-2}, 2N-3,  \dots , 3 ,\frac{1}{2}\] 
	The stabilizer of $\infty$ is generated by $\kappa_{\infty}=\ga_1^{-1}$ that turns in negative direction as $2e^{2\pi i N \lambda}$ does seen as a circle around $\infty$.
	We get an order of cusps around $\infty$:
	\[ 1,0, 2N-1,2N-2, \dots , 3,2\]

	The correspondence is not unique because of symmetries.
	We decide to identify $0 \longleftrightarrow (0:1:1)$ and $1 \longleftrightarrow (1:0:1)$. Then the other correspondences are fixed and we get the claim.
	\end{prf}

\section{Modular forms for Fermat curves}

The modular function and forms for $\G(2)$ are modular function and forms for $\G_N$ as well. Based on $\lambda(z), (1-\lambda)(z)$ 
and $\theta^2(z)$ (see Equation \eqref{eq:prodmodforms}) we can construct further modular function and forms for $\G_N$.

\begin{lemma}
    \label{le:modformgp}
    Let $\lambda(z)$ and $(1-\lambda)(z)$ the modular functions introduced in Equation \eqref{eq:prodmodforms}. The $N$-th roots
    \[ x:=\sqrt[N]{\lambda} \quad y:=\sqrt[N]{1-\lambda} \]
    exist and they are modular functions for $\G_N$.
\end{lemma}
    \begin{prf}
    See \cite{rohrlich}.
    \end{prf}

\begin{lemma}
    \label{le:modformsdiv}
    Remember the cusps $a_j, b_j$ and $ c_j$ of $\G_N$ in Lemma \ref{le:BelyiFermat}. We have the following modular functions and forms for $\G_N$ with divisors as stated.
    \begin{align*}
    \div\, \theta^2 &= \sum_{j=0}^{N-1} N b_j\\
    \div\,x &= \sum_j a_j - \sum_j c_j \\
    \div\,y &= \sum_j b_j - \sum_j c_j \\
    \div\,(x-\zeta^j) &= Nb_j - \sum_j c_j \\
    \div\,(y-\zeta^j) &= Na_j - \sum_j c_j \\
    \div\,(x-\epsilon\zeta^jy) &= Nc_j - \sum_j c_j
    \end{align*} 
    Here we have $\zeta=e^{2\pi i /N}$, $\epsilon=e^{\pi i /N} $,
    $x$ as well as $y$ are from Lemma \ref{le:modformgp} and $\theta^2(z)$ as in Equation \eqref{eq:prodmodforms}. 
\end{lemma}
    \begin{prf}
    See \cite{rohrlich} and \cite{yang}.
    \end{prf}

Now, we can construct modular forms with special zeros.
\begin{lemma}
    \label{le:modformsnn}
    For $j=0, \dots, N-1$ we define
    \begin{align}
      f^{\G_N}_{a_j} &:=\frac{(y-\zeta^j)^N}{y^N}\theta^2 \\
      f^{\G_N}_{b_j} &:=\frac{(x-\zeta^j)^N}{y^N}\theta^2 \\
      f^{\G_N}_{c_j} &:=\frac{(x-\epsilon \zeta^j y)^N}{y^N}\theta^2,
    \end{align}
    where $\zeta=e^{2\pi i /N}$ and $\epsilon=e^{\pi i /N} $. \\
    These all are modular forms for $\G_N$ of weight $2$ and
    \[ \div f^{\G_N}_{i_j} = N^2 i_j, \]
    where $i_j\in\{a_j, b_j, c_j\}$ stands for a cusp of $\G_N$ (see Lemma \ref{le:BelyiFermat}).
\end{lemma}
    \begin{prf}
    Follows easily from Lemma \ref{le:modformsdiv}.
    \end{prf}

By regarding the product of $q$-expansions in several cusps, we recover modular forms for $\G(2)$:

\begin{lemma}
    \label{le:nmformen}
    We build products for the modular forms from Lemma \ref{le:modformsnn} under the action of the slash operator. Thereby, we get only three different results. They are for all $j\in \{0,1, \dots , N-1\}$
    \begin{align}
     \prod_{\g\in \G_N\setminus \G(2)}f^{\G_N}_{a_j}|_{_2}\g(z) & = (-1)^{N^2}\theta^{2N^2}(z)\left(\frac{\lambda(z)}{1-\lambda(z)}\right)^{N^2} \\
     \prod_{\g\in \G_N\setminus \G(2)}f^{\G_N}_{b_j}|_{_2}\g(z) & = (-1)^{N^2}\theta^{2N^2}(z) \\
     \prod_{\g\in \G_N\setminus \G(2)}f^{\G_N}_{c_j}|_{_2}\g(z) & = \theta^{2N^2}(z)\left(\frac{1}{1-\lambda(z)}\right)^{N^2}.
    \end{align}
\end{lemma}
        \begin{prf}
	We have to apply all matrices from \eqref{eq:cosetsG2} to the $f_{i_j}$. 

	The transformational behavior of the form $\theta^2$ is known since $\theta^2\in M_2(\G(2))$. The behavior of $x$ and $y$ is (according to \cite{yang})
	\[  \begin{array}{rclcrcl}
           x|_0\ga_1 &=& \zeta^{-1} x & \quad & x|_0\ga_2 &=& \zeta^{-1}x\\
           y|_0\ga_1 &=& \zeta^{-1}y & \quad & y|_0{\ga_2} &=&  y.
           \end{array} \]
        With this information direct calculations yield the claim.
        \end{prf}

\section{Kronecker limit formulas for Fermat curves}

We want to establish Kronecker limit formulas for the group $\G_N$. From Proposition \ref{prop:KLF} we know that the modular forms introduced in Lemma \ref{le:modformsnn} are suitable.
Missing is the constant $A$ that occurs in Proposition \ref{prop:KLF}.

To calculate that constant, we will use the following trick: The constant $A$ can be calculated by comparing Fourier expansions. 
Since the constant is independent of the cusp in which the Fourier expansion is taken, we may add several expansions and get a multiple of $A$. 
By such a procedure we can obtain $A$ because a suitable sum of expansions leads to known formulas.

\begin{lemma}
    \label{le:limitsum}
    Let $\G_N\subset \G(2)$ be the subgroup associated to the $N$-th Fermat curve, $S_j$ a cusp of $\G_N$ and $f^{\G_N}_j$ the modular form for the cusp $S_j$ according to Lemma \ref{le:modformsnn}.
    It holds:
    \begin{multline}
    \sum_{\g\in \G_N\setminus \G(2)}4\pi \lim_{s\rightarrow 1}\left(E_j^{\G_N}(\g z,s)-\frac{1}{vol(\G_N)(s-1)}\right)
	 = \\ -\frac{1}{N^2}\sum_{\g\in \G_N\setminus \G(2)}\log||f_j^{\G_N}|_{_0}\g(z)||^2 + 4\left( \frac{\zeta'(-1)}{\zeta(-1)} -\log(4\pi) +1+\frac{1}{6}\log(2)-\frac{1}{2}\log(N) \right)
		\label{eq:limitsum}
    \end{multline}
\end{lemma} 
	\begin{prf}
	The left hand side of Equation \eqref{eq:limitsum} had been calculated in Corollary \ref{cor:sumexpgreen}:
	\[4\pi \lim_{s\rightarrow 1}\left(E_j^{\G(2)}(z,s)-\frac{1}{vol(\G(2))(s-1)}\right)-2\log(N) \]

	For the right hand side we realize that \[ \sum_{\g\in \G_N\setminus \G(2)}\log||f_j^{\G_N}|_{_0}\g(z)||^2=\log|| \prod_{\g\in \G_N\setminus \G(2)}f_j^{\G_N}|_{_0}\g(z)||^2.\] Then we use Lemma \ref{le:nmformen} to see that the sum yields $N^2 \log|| G_j(z)||^2$ (the form $G_j$ is one from equations \eqref{eq:mofformG0} to \eqref{eq:mofformGinf}). Therefore we get the claim by comparing the formulas here with the Kronecker limit formula for $\G(2)$ (Equation \eqref{eq:KLFG2}).
	\end{prf}

From the sum formula \eqref{eq:limitsum} we derive individual Kronecker limit formulas.
\begin{thm}[Kronecker limit formula  for $\G_N$]
     \label{thm:grenzeinzeln}
     Let $S_j$ be a cusp of $\Gamma_N$, the subgroup associated to the $N$-th Fermat curve (see Lemma \ref{le:grpfermat}), and let  ${ f^{\G_N}_j\in M_2(\G_N)}$ be the corresponding modular form defined in Lemma \ref{le:modformsnn}.
     We have
     \begin{multline}
        4\pi \lim_{s \rightarrow 1}\left(E_j^{\G_N}(z,s)- \frac{1}{vol(\G_N)(s-1)}\right)= \\
        -\frac{1}{N^2}\log || f^{\G_N}_{j}(z) ||^2
        +\frac{4}{N^2}\left( \frac{\zeta'(-1)}{\zeta(-1)}-\log(4\pi)+1+\frac{1}{6}\log(2)-\frac{1}{2}\log(N)\right).
	\label{eq:grenzeinzeln}
     \end{multline}
\end{thm}
	\begin{prf}
	From Proposition \ref{prop:KLF} follows, that there is an identity
	\[
        4\pi \lim_{s \rightarrow 1}\left(E_j^{\G_N}(z,s)- \frac{1}{vol(\G_N)(s-1)}\right)=
        -\frac{1}{N^2}\log || f^{\G_N}_{j}(z) ||^2+A,
        \]
	where $A$ is a constant. Then Remark \ref{rem:Aindepcusp} explains that $A$ is $\frac{1}{N^2}$ times the constant that occurs in Lemma \ref{le:limitsum}.
	\end{prf}

\section{Scattering constants}

If we take the Fourier expansions in Equation \ref{eq:grenzeinzeln} and compare coefficients, we can get the scattering constants for $\G_N$. 
We only know the Fourier expansion for $f^{\G_N}_j$ in the cusp $\infty$. But because of symmetries of the Fermat curve, this is sufficient to get all scattering constants. 

\begin{thm}
        \label{thm:SCFermatall}
        The scattering constants for $\G_N$, the subgroup associated to the $N$-th Fermat curve, are: \\
        If both cusps are the same, then
        \begin{align}
        C^{\G_N}_{jj}  =\frac{1}{6N^2}\left(C^{\Gamma(1)}- \frac{1}{\pi}\left((12N+2)\log(2)+(-3N+6)\log(N)\right)\right). \label{eq:SCfermatjj}
        \end{align}
        If $S_j\neq S_k$ and $\belyi_N(S_j) \neq \belyi_N(S_k)$, then
        \begin{align}
        C^{\G_N}_{kj} =\frac{1}{6N^2}\left( C^{\Gamma(1)}- \frac{1}{\pi}\left(2\log(2)+6\log(N) \right) \right) .
             \label{eq:SCfermatkj}
        \end{align}
        If $S_j\neq S_l$ but $\belyi_N(S_j) = \belyi_N(S_l)$, then
        \begin{align}
        C^{\G_N}_{lj} =\frac{1}{6N^2} \left( C^{\Gamma(1)}- \frac{1}{\pi}\left(2\log(2)+6\log(N)
            +3N\log\big|1-\zeta_N^{l-j} \big|\right)\right) .
            \label{eq:SCfermatlj}
        \end{align}
        The number $\zeta_N^{l-j}=\zeta_N^{l}\left(\zeta_N^{j}\right)^{-1}$ is given by the $N$-th roots of unity which determine the cusps $S_l$ and $S_j$ in Lemma \ref{le:BelyiFermat}.
\end{thm}
	\begin{prf}
	In the case that the second cusp $S_j=\infty$ we get the scattering constants via comparison of the first coefficients in the
	realizations of Equation \eqref{eq:grenzeinzeln}. We expand both sides of \eqref{eq:grenzeinzeln} in $\infty$ and calculate the constant term:

	The constant term of the Eisenstein series at $s=1$ is (see Equation \eqref{eq:expgreen})
	\[ \delta_{jk}\frac{2\pi y}{N}+4\pi \widetilde{C}_{jk}^{\G_N}-\frac{2}{N^2}\log(y).\]
	On the other side of Equation \eqref{eq:grenzeinzeln} we need the constant term of the $q$-expansion of $f_*^{\G_N}$. 
	We get the constant term if we look at the $q$-expansions of $\theta^2$, $\lambda$ as well as $1-\lambda$ and calculate their roots. The result is
	\[ \text{constant term}\left(\log||f_j^{\G_N}||^2\right)=\begin{cases}
		2\log(y) -2N\pi y +N(\log(16)-\log(N)) & \text{ if } S_j=c_0 \\
		2\log(y)+\log|1-\zeta_N^j| & \text{ if } S_j=c_j  \\ 
		2\log(y) & \text{ if } S_j\in\{a_j,b_j\},
	\end{cases}\]
	where the names of the cusps come from Equation \eqref{eq:cuspsfermat}. 
	The root $\zeta^j$ in the second case is determined by the cusp $c_j= (\epsilon \zeta_N^j:1:0)$ ($j\in \{1, \dots , N-1\}$). 
		
	By taking the constant from Theorem \ref{thm:grenzeinzeln} and the difference $\widetilde{C}_{jk}^{\G_N}= C_{jk}^{\G_N} +\frac{1}{2N^2\pi}\log(2N)$ into account, we derive all formulas from the statement.

	Now, all that remains to be chown is that these formulas generalize to arbitrary second cusp. 

	Via the definition of Eisenstein series it is easy to show that for $\g\in\G(1)$
	\[ E_{\g^{-1}(j)}^{\G}(\g_l^{-1}z,s)=E_j^{\g^{-1}\G\g}(\g_l^{-1}\g^{-1}z,s)\] 
	holds. Since $\G_N\triangleleft \G(1)$, we have $\G_N=\g^{-1}\G_N\g$ and derive
	\[ C^{\G_N}_{jl}= C^{\G_N}_{\g(j),\g(l)}.\]

	Therefore, to get all scattering constants we have to find elements of $\G(1)$ that map $\infty$ to all cusps
	and to find out to where these matrices move the second cusp. 

	The conjugation with $S=\mat{0}{1}{-1}{0}$ (to get $0$)  or $TS=\mat{1}{1}{0}{1}\mat{0}{1}{-1}{0}$ (to get $1$) give (for even $j,l$ and odd $k$ in the range $1,2,\dots ,2N$)
		\begin{align*}
		C^{\G_N}_{1/l, \infty} &= C^{\G_N}_{2N-l, 0} &	C^{\G_N}_{1/l, \infty} &= C^{\G_N}_{2N-l+1, 1} \\
		C^{\G_N}_{j, \infty} &= C^{\G_N}_{1/(2N-j), 0} & C^{\G_N}_{j, \infty} &= C^{\G_N}_{1/j, 1} \\
		C^{\G_N}_{k, \infty} &= C^{\G_N}_{2N-k, 0} 	&  C^{\G_N}_{k, \infty} &= C^{\G_N}_{2N-k+1, 1}. 
		\end{align*}
        (To see this we use 
		\begin{align*}
		-l &\sim_{_{\G_N}} 2N-l &\qquad \text{ via }&\qquad \ga_1^N \\
		-1/j&\sim_{_{\G_N}} 1/(2N-j)  &\qquad \text{ via }&\qquad \ga_2^N \\
		-1/k&\sim_{_{\G_N}} 2N-k  & \qquad \text{ via }&\qquad \ga_1^{(2N+k-1)/2}(\ga_2\ga_1^{-1})^{(k-1)/2}\ga_2^{(1-k)/2} \\
		(j-1)/j&\sim_{_{\G_N}} 1/j & \qquad \text{ via }&\qquad \ga_2^{j/2}\ga_1^{(j-2N)/2}(\ga_1\ga_2^{-1})^{j/2} \\ 
		(k-1)/k&\sim_{_{\G_N}} 2N-k+1 & \qquad \text{ via }&\qquad \ga_1^{(2N-k+1)/2}\ga_2^{(k-1)/2}(\ga_2\ga_1^{-1})^{(1-k)/2}
		\end{align*}
        to identify the cusps.)

	Hence, all scattering constants for $0,1$ or  $\infty$ as second cusp are known.  

	With $\g_1$ in the cases $0$ as well as $1$ and $\g_2$ in the case of $\infty$ we generalize to all cusps. (In the only cases we really need, i.e. the ones where both cusps involved are equivalent under $\G(2)$, 
	there occur no further difficulties when we try to identify the resulting cusps in the system of representatives.)  

	Thereby, we will find the formulas from the statement when we use the correspondence of cusps in 
	Proposition \ref{prop:corrcuspsfermat} and realize that $|1-\zeta_N^{j}|=|1-\zeta_N^{N-j}|$. 
	\end{prf}

\begin{rem}
    There is an alternative method to determine the scattering constants for the Fermat curves by means of Arakelov theory. 
    U. K\"uhn \cite{kuehn} showed how scattering constants give arithmetic intersection numbers in the infinite places. 
    Together with the intersection numbers in the finite places that Ch.~Curilla \cite{curilla} calculated in his PhD-thesis, we could get the scattering constants without using Kronecker limit formulas and the symmetries of the Fermat curves by studying arithmetic intersection numbers in the cusps.
\end{rem}

\begin{rem}
	We can approximate scattering constants for the Fermat curves numerically. 
	The means for that has been developed by the author in her Diplomarbeit \cite{posingiesdipl} and her Dissertation \cite{posingiesdiss}.
        For some small $N$, the structure of the scattering matrix, i.e.~the symmetries, and the values of the scattering constants were replicated numerically.
\end{rem}

\end{document}